\documentclass[12pt]{amsart}

\usepackage{srcltx}
\usepackage[varg]{txfonts}
\usepackage[all]{xy}
\usepackage{amsfonts, amsmath, amssymb, graphics, graphicx, caption}

\newtheorem{thm}{Theorem}[section]
 
 \newtheorem{cor}[thm]{Corollary}
 \newtheorem{lem}[thm]{Lemma}
 \newtheorem{prop}[thm]{Proposition}
 
 \newtheorem{dfn}[thm]{Definition}

 \newtheorem{rmk}[thm]{Remark}
 \newtheorem{ex}[thm]{Example}
 \theoremstyle{definition}
 \theoremstyle{remark}
 \numberwithin{equation}{section}

\newcommand\CC{\mathbb{C}}
\newcommand\NN{\mathbb{N}}
\newcommand\QQ{\mathbb{Q}}
\newcommand\RR{\mathbb{R}}
\newcommand\ZZ{\mathbb{Z}}
\newcommand\HH{\mathbb{H}}

\newcommand\FF{\mathfrak{F}}
\newcommand\reg{\mathrm{reg}}

\newcommand{\sm}{\left(\begin{smallmatrix}}
\newcommand{\esm}{\end{smallmatrix}\right)}

\newcommand\re{\mathrm{Re}\,}
\newcommand\im{\mathrm{Im}\,}
\newcommand\SL{{\mathrm{SL}}}

\newcommand\poiss[1]{\mathsf{P}_{\!s}}
\newcommand\tpoiss[1]{\mathsf{P}^\dagger_{\!s}}

\makeatletter
\newcommand\hmatc[4]{\left[ {#1\@@atop #3}{#2\@@atop #4}\right]}
\newcommand\hmatr[4]{\left[ {\hfill #1\@@atop\hfill #3}{\hfill
#2\@@atop\hfill #4}\right]}
\newcommand\matc[4]{\left( {#1\@@atop #3}{#2\@@atop #4}\right)}
\newcommand\matr[4]{\left( {\hfill #1\@@atop\hfill #3}{\hfill
#2\@@atop\hfill #4}\right)}
\makeatother

\newcommand\be{\begin{equation}}
\newcommand\ee{\end{equation}}
\newcommand\bad{\be\begin{aligned}}
\newcommand\ead{\end{aligned}\ee}
\newcommand\badl[1]{\be\label{#1}\begin{aligned}}

\begin{document}
\title{Theta liftings on weak Maass forms}
\author{YoungJu  Choie }

 \address{Department of Mathematics and PMI\\
 Pohang University of Science and Technology\\
 Pohang, 790--784, Korea}
 \email{yjchoie@gmail.com}

 \author{Subong Lim }

 \address{Department of Mathematics and PMI\\
 Pohang University of Science and Technology\\
 Pohang, 790--784, Korea}
  \email{subong@postech.ac.kr}

 \thanks{Keynote: theta lifting, weak Maass form }


\begin{abstract}
We construct theta liftings from half-integral weight weak Maass
forms to even integral weight weak Maass forms by using regularized theta integral. Moreover it gives
an extension of Niwa's theta liftings on harmonic weak Maass
forms. And we obtain the similar results to those by Niwa.
\end{abstract}

\maketitle

\raggedbottom

\section{\bf{Introduction}}

 Weak Maass forms are Maass forms which allow exponential growth at cusps.
 We denote the space of weak Maass forms of weight $k$ and eigenvalue
  $s$ for a congruence subgroup $\Gamma$ and a character $\chi$ as
  $WMF_{k,s}(\Gamma,\chi)$. And harmonic weak Maass forms are weak Maass forms
  whose eigenvalue is zero. We denote the space of harmonic weak Maass forms of
   weight $k$ for a congruence subgroup $\Gamma$ and a character $\chi$ as
    $H_{k}(\Gamma,\chi)$. Harmonic weak Maass forms are related to Ramanujan's
    mock theta functions. A mock theta function was introduced by Ramanujan in
    his last letter to Hardy. He gave some examples and definition. In 2001,
    Zwegers\cite{Z2} discovered that mock theta functions are holomorphic parts
    of nonholomorphic modular forms of weight $1/2$. More generally, Zagier\cite{Z1}
    defined mock modular forms. And it turns out that mock modular forms
    are holomorphic parts of harmonic weak Maass forms of any
    weight.\\

 In this paper, we construct even-integral weight weak Maass forms from half-integral
 weight weak Maass forms. To construct maps, we shall use the indefinite theta series.
 This theta series was used by Shintani, Niwa and Cipra. Using this theta series, Shintani\cite{S2}
 constructed an inverse map of Shimura correspondence and Niwa\cite{N} gave another method
 to construct even-integral weight cusp forms from half-integral weight cusp forms and
 proved the conjecture about the level in the Shimura correspondence. And Cipra\cite{C}
 extended Niwa's lifting to all holomorphic modular forms of all positive, half-integral weight.
 And we will regularize theta integral to make theta liftings well-defined in the case of weak Maass forms. Borcherds\cite{B1}
 used a regularized theta integral and constructed theta liftings on weakly holomorphic modular
 forms. And Bruinier\cite{B3} defined theta lifting on harmonic weak Maass forms of weight $1/2$.
 \\

 Let $\HH$ be the upper half-plane and let $z = u+iv, w = \xi+i\eta\in \HH$.
 Following Cipra's method in \cite{C} for given $N\in\NN$ and a character $\chi$ for $\Gamma_0(4N)$ we obtain theta functions $\theta(z,w;f_{k,m})$
  where $k,m\in\ZZ$. This theta function has an important property: $\theta(z,w;f_{k,m})$ is a nonholomorphic
  modular form of weight $k/2$ for $\Gamma_0(4N)$ and $\chi$ as a function of $z$ and $\overline{\theta(z,w;f_{k,m})}$
  is a nonholomorphic modular form of weight $2m$ for $\Gamma_0(2N)$ and $\chi^2$ as a function of $w$. For a given weak Maass form $g$ of weight $k/2$ for $\Gamma_0(4N)$ and $\chi$ we define a function on $\HH$
 \[\Phi(g)(w) = \int^{\reg}_{\SL_2(\ZZ)\setminus\HH}v^{k/2}
 \sum_{\alpha\in\Gamma_0(4N)\setminus\SL_2(\ZZ)}g_\alpha(z)\overline{\theta_\alpha(z,w;f_{k,m})}\frac{dudv}{v^2}\]
 and
 \[\Phi_D(g)(w) = \Phi(g_D)(\frac wD)\]
 where $(g_\alpha)(z) = (cz+d)^{-k/2}g(\frac{az+b}{cz+d})$ for $\alpha = \sm a&b\\c&d\esm$ and $g_D(z) = g(Dz)$.
 Here $\int^{\mathrm{reg}}_{\SL_2(\ZZ)\setminus\HH}$ is a regularized integral which
 was used by Borcherds \cite{B1}.
 Then this gives a lifting from weak Maass forms of weight $k/2$ to weak Maass forms
 of weight $2m$. The following is the precise statement:

 \begin{thm} \label{main1} Let $k$ and $m$ be any integers and let $g(z)\in WMF^*_{k/2,s}(\Gamma_0(4N),\chi)$.
 Assume that $\chi(0)=0$. Then $\Phi_D(g)(w)$ is a weak Maass form of weight $2m$ and
 eigenvalue $4m(m-1)-3-k(k-4)+4s$ for $\Gamma_0(2N,D)$ and $\chi^2$. But
 it may have singularities at Heegner points of the form $w = \frac{b\pm\sqrt{b^2-ac}}{2Nc}\in\HH$ for
 $a,b,c\in\ZZ$, while $\Phi_D(g)$ does not have any singularities
 in the case when $g$ is a usual Maass form and $k\neq1$.
\end{thm}

\begin{rmk} Shimura correspondence is the map from cusp
 forms of weight $k/2$ to cusp forms of weight $2\lambda=k-1$.
 Unlike the Shimura correspondence, we can choose $k$ and $m$ independently. And $k$ and $m$ have any sign.
 \end{rmk}

Here $ WMF^*_{k/2,s}(\Gamma_0(4N),\chi) $  is a subspace of $
WMF_{k/2,s}(\Gamma_0(4N),\chi) ,$
 the space of weak Maass forms, whose precise definition will be given in Section 2.
 And $\Gamma_0(2N,D) = \{\sm a&b\\c&d\esm\in \Gamma_0(2N)|\ b\equiv0
 (D)\}$. In particular, Harmonic weak Maass forms are weak Maass forms
whose eigenvalues are zero. The space of
  harmonic weak Maass forms are denoted by
 $ H_{k/2}(\Gamma_0(4N),\chi)$, that is,
  $H_{k/2}(\Gamma_0(4N),\chi)=WMF_{k/2,s}(\Gamma_0(4N),\chi)$.
 And we have subspaces $H^*_{k/2}(\Gamma_0(4N),\chi)$ and $ H^+_{k/2}(\Gamma_0(4N),\chi)$ of $H_{k/2}(\Gamma_0(4N),\chi)$,
 which will be defined in Section 2. Note that if $g(z)\in H^+_{k/2}(\Gamma_0(4N),\chi)$
  then $g$ has a Fourier expansion of the form
 \[g(z) = \sum_{n \in \mathbb{Q} \atop n\gg -\infty}a^+(n)e(nz)+ \sum_{n \in \mathbb{Q} \atop n<0}a^-(n)W_{k/2}(2\pi nv)e(nz),\]
 where $W_k(x) := \int^\infty_{-2x}e^{-t}t^{-k}dt = \Gamma(1-k,2|x|)$ for $x<0$ and $e(z) = e^{2\pi iz}$.

 \begin{thm} \label{main2} Let $k\geq1$ and $\lambda = \frac{k-1}{2}$. Assume $\chi(0)=0$ and
 that $a^+(0) = 0$ if $k=1$ and $\chi$ is a principal character. And
 let $g\in H^*_{k/2}(\Gamma_0(4N),\chi)$. If $g(z)$ has a Fourier expansion as above then
 $\Phi_D(g)(w) \in H^*_{2\lambda}(\Gamma_0(2N,D),\chi^2)$ with the same possible singularities at Heegner points
 as in Theorem \ref{main1}. Moreover if $g\in H^+_{k/2}(\Gamma_0(4N),\chi)$ then we have the following results:

 \begin{enumerate}
 \item If $k\geq3$ then we have
 \[\Phi_D(g)(i\infty) = \lim_{\eta\to\infty}(\Phi_D(g)(i\eta)) = C_D(\lambda)\frac{a^+(0)}{2}L(1-\lambda,\chi_D)\]
 and if $k=1$ then we have
 \[\Phi(g)(i\infty) = 4N^{1/4}a^+(0)\sum_{m=1}^\infty \frac{\chi(m)}{m}.\]

 \item If $k\geq3$ then we have
 \[\int^\infty_0 \eta^{s-1}(\Phi_D(g)(i\eta)-\Phi_D(g)(i\infty))d\eta\]
 \[=  C_D(\lambda)(2\pi)^{-s}\Gamma(s)L(s+1-\lambda,\chi_D)\sum_{n=1}^\infty\frac{a^+(Dn^2)}{n^s}\]
 and if $k=1$ then we have
 \[\int^\infty_0 \eta^{s-1}(\Phi(g)(i\eta)-\Phi(g)(i\infty))d\eta =
 C_1(\lambda)(2\pi)^{-s}\Gamma(s)L(s+1,\chi)\sum_{n=1}^\infty\frac{a^+(n^2)}{n^s}\]
 \end{enumerate}
 where $C_D(\lambda) = (-1)^\lambda 2^{-3\lambda+2}(DN)^{\lambda/2+1/4}$
 and $\chi_D = \chi(\frac{-1}{})^\lambda(\frac{D}{})$.
 \end{thm}

  Then this is the analogous result with Niwa\cite{N} and Cipra\cite{C}. So this is
  the extension of Niwa's theta lifting to harmonic weak Maass forms.

 \begin{rmk} If $g$ is a usual Maass form and $k\neq1$,  then we do not need to regularize integral because
 $\theta(z,w;f_{k,\lambda})$ is rapidly decreasing at all cusps (see Theorem 2.6
 in \cite{C}).
 \end{rmk}

 \begin{ex}

 We consider a holomorphic Eisenstein series (see Proposition 2.10 in \cite{C}): Let $k\geq3$. Define
 \[E_{k/2}(z,s) = E_{k/2}(z,s,4N,\chi) = \sum_{\gamma\in \Gamma_\infty\setminus\Gamma_0(4N)}\bar{\chi}(d)\frac{\im(\gamma z)^s}{j(\gamma,z)^k}.\]
 Then it is a Maass form of weight $k/2$ with eigenvalue $-s(s-1)-ks/2$ for $\Gamma_0(4N)$ and $\chi$. Assume that $\chi_1 = \chi(\frac{-1}{})^\lambda$ is primitive mod $4N$. Recall that $\lambda = (k-1)/2$. Then
 \[\Phi_1(E_{k/2}(z,s))(w) = C(s)E_{2\lambda}(w,2s)\]
 where
 \[E_{2\lambda}(w,2s) = E_{2\lambda}(w,2s,2N,\chi^2) = \sum_{\gamma\in\Gamma_\infty\setminus\Gamma_0(2N)}\bar{\chi}^2(d)\frac{\im(\gamma z)^{2s}}{(cw+d)^{2\lambda}}.\]
 Then it is a Maass form of weight $2\lambda$ with eigenvalue $-2s(2s-1)-4\lambda s$ for $\Gamma_0(2N)$ and $\chi^2$.

 \end{ex}

 \begin{thm} Let $g$ be the same as that in Theorem \ref{main2}.
 If $\Phi_D(g)(w) \in H^+_{2\lambda}(\Gamma_0(2N,D),\chi^2)$ and we write its Fourier expansion as follows
 \[\Phi_D(g)(w) = \sum_{n\geq0}A_D^+(n)e(nw) + \sum_{n<0}A_D^-(n)W(2\pi n\eta)e(nw)\]
 then we have
 \begin{eqnarray*}
 &&\sum_{n=1}^\infty \frac{A_D^+(n)}{n^s} + \frac{1}{\Gamma(s)}\big(\int_0^\infty\int_{2y}^\infty
 e^{y-x}x^{-k/2}y^{s-1}dxdy\big)\sum_{n=1}^\infty\frac{A_D^-(-n)}{n^s} \\
 &=& C_D(\lambda)L(s-\lambda+1,\chi_D)\sum_{n=1}^\infty\frac{a^+(Dn^2)}{n^s}.
 \end{eqnarray*}

 \end{thm}
So we have the similar result with Shimura's correspondence in
\cite{S1} but we have extra terms.


 \section{\bf{Weak Maass forms}}


 In this section, we review weak Maass forms and harmonic weak Maass forms, which were
 introduced in \cite{B4} and \cite{B5}. For $\gamma = \sm a&b\\c&d\esm\in \SL_2(\ZZ)$ and $z\in \mathbb{H}$, let
 \[J(\gamma,z) = cz + d,\ j(\gamma, z) = \epsilon_d^{-1}(\frac cd)J(\gamma,z)^{1/2}\ \text{and}\ \gamma z = \frac{az+b}{cz+d},\]
 where $\epsilon_d = 1$ or $i$ as $d\equiv 1$ or $3$, and $(\frac cd)$ is the quadratic residue symbol as defined in \cite{S1}.

 If $k\in\mathbb{Z}, \gamma = \sm a&b\\c&d\esm \in \SL_2(\mathbb{Z})$, and $f$ is a function on $\mathbb{H}$, define
 \[(f|_k\gamma)(z) = (cz+d)^{-k}f(\gamma z).\]
 When $k\in\mathbb{Z}$ is odd and $\gamma = \sm a&b\\c&d\esm \in \Gamma_0(4)$, define
 \[(f|_{k/2}\gamma)(z) = j(\gamma,z)^{-k}f(\gamma z).\]

 Let $f$ be a smooth function on $\mathbb{H}$ and let $k$ be an integer or half-integer.
 And let $\Gamma$ be a congruence subgroup of $\SL_2(\ZZ)$ and let $\chi$ be a character for $\Gamma$. And let
 \[{\bf{G_k(\Gamma,\chi)}} = \{ f\in C^\infty(\HH) |\ f|_k\gamma = \chi(d)f\ \text{for all}\ \gamma = \sm a&b\\c&d\esm \in \Gamma\}.\]
 If $f\in G_k(\Gamma,\chi)$ then we say that $f$ is a nonholomorphic modular form of weight $k$ for $\Gamma$ and $\chi$.

 Define $\Delta_k = -v^2(\frac{\partial^2}{\partial u^2}+\frac{\partial^2}{\partial v^2})
 +ikv(\frac{\partial}{\partial u}+i\frac{\partial}{\partial v})$. This is the usual hyperbolic Laplace operator in weight
 $k$.\\

 \begin{dfn} A twice continuously differentiable function on $\HH$ is called a Maass form
 of weight $k$ and eigenvalue $s$ for $\Gamma$ and $\chi$ if
 \begin{enumerate}
 \item $(f|_k\gamma) = \chi(\gamma)f$ for all $\gamma\in \Gamma$,
 \item $\Delta_k f =s f$ for $s\in\CC$,
 \item For all $\gamma\in \SL_2(\ZZ)$, $(f|_k\gamma)(z) = O(v^\delta)$ as $v\to\infty$ for some $\delta>0$.
 \end{enumerate}
 We denote the vector space of those {\bf{Maass forms}} by ${\bf{Maass_{k,s}(\Gamma, \chi)}}$.
 If we change the third condition, which is a growth condition, into
 \begin{enumerate}
 \item[(3)'] For all $\gamma\in \SL_2(\ZZ)$, $(f|_k\gamma)(z) = O(e^{\delta v})$ as $v\to\infty$ for some $\delta>0$.
 \end{enumerate}
 then $f$ is called a {\bf{weak Maass form}}.  Especially, if $s=0$, then it is called a {\bf{harmonic weak Maass form}}.
 \\

 Let ${\bf{WMF_{k,s}(\Gamma,\chi), H_{k}(\Gamma,\chi)}}$ denote respectively the spaces of weak Maass forms and harmonic weak Maass forms.
 \end{dfn}

 If $f$ is a weak Maass form for $\Gamma$, then it satisfies $f(z+l) = f(z)$ since $\sm 1&l\\0&1\esm\in\Gamma$
 for some positive integer $l$. Hence there is a Fourier expansion
 \begin{equation} \label{fourier0}
 f(z) = \sum_{n\in\QQ}a(n,v)e(nu).
 \end{equation}
 Then we define
 \[{\bf{WMF^*_{k,s}(\Gamma,\chi)}} = \{f\in WMF_{k,s}(\Gamma,\chi)|\ a(n,v)e^{2\pi nv} = O(v^\delta)\
  \text{for some $\delta>0\ \forall n$}\}.\]
 This is the subspace of $WMF_{k,s}(\Gamma,\chi)$.

 In particular, ${\bf{H^*_{k}(\Gamma,\chi)}}$ can be characterized by the differential operator $\xi_k$,
 which is studied by Bruinier and Funke\cite{B4}: We define a differential operator
 $\xi_k = 2iv^k\overline{\frac{\partial}{\partial\bar{z}}}$. This gives a map
 \begin{eqnarray*}
 \xi_k: H_{k}(\Gamma,\chi) &\twoheadrightarrow& M^!_{2-k}(\Gamma,\bar{\chi}),
 \end{eqnarray*}
 where $M^!_k(\Gamma,\chi)$ is the space of weakly holomorphic modular forms of
  weight $k$ for $\Gamma$ and $\chi$. It is easy to see that $H^*_{k}(\Gamma,\chi) = \xi_k^{-1}(M_{2-k}(\Gamma,\bar{\chi}))$.
   Here $M_k(\Gamma,\chi)$ is the space of holomorphic modular forms of weight $k$ for $\Gamma$ and $\chi$.

 The Fourier expansion of any $f \in H^*_{k}(\Gamma,\chi)$ gives a unique decomposition $f = f^+ + f^-$, where
 \begin{eqnarray} \label{fourier}
 f^+(z) &=& \sum_{n \in \mathbb{Q} \atop n\gg -\infty}a^+(n)e(nz),\\
 \nonumber f^-(z) &=& \sum_{n \in \mathbb{Q} \atop n<0}a^-(n)W_k(2\pi nv)e(nz) + a^-(0)v^{1-k/2},
 \end{eqnarray}
 and $W_k(x) := \int^\infty_{-2x}e^{-t}t^{-k}dt = \Gamma(1-k,2|x|)$ for $x<0$ and $e(z) = e^{2\pi iz}$. Here
 $f^+$ and $f^-$ are called the holomorphic part and nonholomorphic part of $f$,
  respectively. Note that $f$ has a Fourier expansion of the form (\ref{fourier}) at all cusps. And the Fourier polynomial
 \[P(f) = \sum_{n\leq0}a^+(n)q^n\ \text{with}\ q = e^{2\pi iz}\]
 is called the principal part of $f$.

 We define another subspace ${\bf{H^+_{k}(\Gamma,\chi)}}$ of $H^*_{k}(\Gamma,\chi)$ as
 the inverse image of $S_{2-k}(\Gamma,\bar{\chi})$ under the map $\xi_k$. So it has a Fourier expansion
 \begin{equation} \label{fourier+}
 f(z) = \sum_{n \in \mathbb{Q} \atop n\gg -\infty}a^+(n)e(nz) + \sum_{n \in \mathbb{Q} \atop n<0}a^-(n)W_k(2\pi nv)e(nz).
 \end{equation}

\medskip

In summary we have introduced the following notations:
$$WMF_{k,s}^*(\Gamma, \chi) \subset WMF_{k,s}(\Gamma, \chi)
\subset Maass_{k,s}(\Gamma, \chi) \subset G_{k}(\Gamma, \chi)  $$

$$H^+_{k}(\Gamma,\chi)\subset H^*_{k}(\Gamma,\chi)\subset H_{k}(\Gamma, \chi)=  WMF_{k,0}(\Gamma,
\chi)$$

\bigskip

 \section{{\bf{Construction of Theta lifting}}}

 \subsection{Indefinite theta functions}

 In this section we define theta functions by using Cipra's method in \cite{C}.
  The most of results in this section are in \cite{C}. We begin by defining
  the Weil representations on $\SL_2(\mathbb{R})$ on the space of Schwartz
  functions $\mathcal{S}(\mathbb{R}^n)$. Let $Q$ be a rational symmetric
  matrix of signature $(p,q),p+q=n$. For $x,y\in \mathbb{R}^n$, define the inner product
 \[<x,y>\ =\ ^txQy.\]
 For a matrix $\gamma=\sm a&b\\c&d\esm\in \SL_2(\mathbb{R})$ and a Schwartz
 function $f\in \mathcal{S}(\mathbb{R}^n)$, define the Weil representation
 \begin{equation*}
 (r(\gamma,Q)f)(x)=
 \begin{cases}
 |a|^{n/2}e[\frac{ab}2<x,x>]f(ax) & \text{if $c=0$},\\
 |\det Q|^{-1/2}|c|^{-n/2}\int_{\mathbb{R}^n}e[\frac{a<x,x>-2<x,y>+d<y,y>}{2c}]f(y)dy & \text{if $c\neq0$}.
 \end{cases}
 \end{equation*}

 \begin{prop}\label{matrix} Let $\gamma\in \SL_2(\mathbb{R})$ and $\sigma_z =
  \sm v^{1/2}& uv^{-1/2}\\ 0& v^{-1/2}\esm$ for $z = u+iv\in\mathbb{H}$.
  Define $\phi (\mathrm{mod}\ 2\pi)$ by $e^{-i\phi} = J(\gamma,z)/|J(\gamma,z)|$, and let
 \[\kappa(\phi) = \sm \cos\phi&\sin\phi\\-\sin\phi&\cos\phi\esm.\]
 Then
 \begin{enumerate}
 \item $\gamma\sigma_z = \sigma_{\gamma z}\kappa(\phi)$,
 \item $r(\gamma,Q)r(\sigma_z,Q) = r(\sigma_{\gamma z},Q)r(\kappa(\phi),Q)$.
 \end{enumerate}
 \end{prop}

 {\bf{Proof}} This is Proposition 1.3 of \cite{C}. \qed

 \begin{cor} \label{second} For $\gamma = \sm a&b\\c&d\esm\in \SL_2(\mathbb{R})$, and $t\in\mathbb{R}$,
 let $\gamma_t = \sm a&bt^2\\c/t^2&d\esm = \sm t & \\ & t^{-1}\esm \gamma \sm t^{-1}&\\&t\esm$. Let $\kappa(\phi)$ be as before. Then
 \begin{equation} \label{p}
 \gamma_t\sigma_{t^2z} = \sigma_{t^2(\gamma z)}\kappa(\phi)
 \end{equation}
 and
 \[r(\gamma_t,Q)r(\sigma_{t^2z},Q) = r(\sigma_{t^2(\gamma z)},Q)r(\kappa(\phi),Q).\]
 \end{cor}

 {\bf{Proof}} This follows since $J(\gamma_t,t^2z)/|J(\gamma_t,t^2z)| = J(\gamma,z)/|J(\gamma,z)|$. \qed

 \medskip

 Let $L$ be an even lattice in $\mathbb{R}^n$ and let $L^*$ be the dual lattice.
  Denote by $v(L)$ the volume of a fundamental parallelotope of $L$ in $\mathbb{R}^n$:
 \[v(L) = \int_{\mathbb{R}^n/L}dx.\]
 Let $\{\lambda_1,\cdots,\lambda_n\}$ be a $\ZZ$-basis for $L$, and define $B = \det(<\lambda_i,\lambda_j>)$.

 \begin{dfn}
 \begin{enumerate}
 \item We say that a function $\omega: L^*/L\to\CC$ has the {\bf{first permutation property}}
 for $\Gamma_0(4N)$ with a character $\chi$ if it satisfies
 \begin{enumerate}
 \item $\omega(\kappa) =0$ if $<\kappa,\kappa>\not\in 2\mathbb{Z}$,
 \item $\omega(d\kappa) = \chi(d)\omega(\kappa)$ for $\gamma = \sm a&b\\c&d\esm \in \Gamma_0(4N)$,
 \end{enumerate}
 where $\chi$ is a character mod $4N$.

 \item We say that a function $f\in \mathcal{S}(\mathbb{R}^n)$ has the {\bf{first spherical property}}
 for weight $k/2$ if it satisfies
 \[r(\kappa(\phi),Q)f = \epsilon(\kappa(\phi))^{p-q}\sqrt{e^{-i\phi}}^{-k}f,\]
 for all $\kappa(\phi) = \sm \cos\phi& \sin\phi \\ -\sin\phi & \cos\phi \esm$ where $k\in\mathbb{Z}$ and for $\gamma = \sm a&b\\c&d\esm$
 \begin{equation*}
 \epsilon(\gamma) =
 \begin{cases}
 \sqrt{i} & c>0,\\
 i^{(1-\operatorname{sgn} d)/2} & c = 0,\\
 \sqrt{i}^{-1} & c<0.
 \end{cases}
 \end{equation*}

 \end{enumerate}
 \end{dfn}

 Let $f\in\mathcal{S}(\RR^n)$ and define, for $h\in L^*/L$
 \[\theta(f,h) := \sum_{x\in L}f(x+h).\]
 Then , by Shintani\cite{S2}:

 \begin{prop} \label{Shintani0} Let $\gamma = \sm a&b\\c&d\esm\in \SL_2(\ZZ)$. Then
 \[\theta(r(\gamma,Q)f,h) = \sum_{j\in L^*/L}c(h,j)_\gamma\theta(f,j)\]
 where
 \begin{equation*}
 c(h,j)_\gamma =
 \begin{cases}
 \delta_{h,aj}e(\frac{ab}2<h,j>)& \text{if $c=0$},\\
 |\det Q|^{1/2}v(L)^{-1}|c|^{-n/2}\displaystyle\sum_{r\in L/cL}e[\frac1{2c}(a<h+r,h+r>\\
 -2<j,h+r>  + d<j,j>] & \text{if $c\neq0$}.
 \end{cases}
 \end{equation*}
 \end{prop}

 Take $f$ having the first spherical property for weight $k/2$, and
 let $\omega$ have the first permutation property for $\Gamma_0(4N)$ with character $\chi$. Define
 \[\theta(z,f,h) := v^{-k/4}\theta(r(\sigma_z,Q)f,h), \ h\in L^*/L\]
 and
 \[\theta(z,f; \omega) := \sum_{h\in L^*/L}\omega(h)\theta(z,f,h).\]

 \begin{thm} [Shintani] \label{Shintani} Let $\gamma =
 \sm a&b\\c&d\esm\in \SL_2(\mathbb{Z})$. Then
 \[(cz+d)^{-k/2}\theta(\gamma z,f,h) = \sqrt{i}^{-(p-q)
 \operatorname{sgn} c}\sum_{j\in L^*/L}c(h,j)_\gamma\theta(z,f,j)\]
 where $c(h,j)_\gamma$ as in Proposition \ref{Shintani0}.
 \end{thm}

 {\bf{Proof}} This is Theorem 1.5 in \cite{C}. \qed

 \begin{cor} [Shintani] \label{first} Let $\gamma = \sm a&b\\c&d\esm\in\Gamma_0(4N)$. Then
 \[j(\gamma,z)^{-k}\theta(\gamma z, f; \omega) = \chi'(d)\theta(z,f; \omega)\]
 where $\chi'(d) = (\frac{-1}d)^{(k-n)/2}(\frac 2d)^n(\frac Bd)((-1)^qB,d)_\infty\chi(d)$ with the Hilbert symbol
 \begin{equation*}
 (x,y)_\infty =
 \begin{cases}
 -1 & \text{if $x,y<0$},\\
 1 & \text{otherwise}.
 \end{cases}
 \end{equation*}
 \end{cor}

 \medskip

 Let $O(Q)$ be the orthogonal group of $Q$: $O(Q) = \{g|\ ^tgQg = Q\}$. Let $SO(Q)$ denote the
 connected component of the identity in $O(Q)$, consisting of those matrices $g$
 with $\det g=1$. We define a unitary representation of $SO(Q)$ on
 $L^2(\mathbb{R}^n)$ by letting $(p(g)f)(x) = f(g^{-1}x)$.
 By definition of $SO(Q)$, $p(g)$ commutes with the Weil representation (See \cite{C} page 64):
 \begin{equation} \label{commute}
 p(g)(r(\gamma,Q)f) = r(\gamma,Q)(p(g)f).
 \end{equation}

 Now we introduce theta kernel as in \cite{C}. Take the following special Q: let
 \begin{equation} \label{Q}
 Q = \frac2N\sm &&-1\\&1&\\-2&&\esm,
 \end{equation}
 be a matrix with signature $(2,1)$. Let $L = 4N\ZZ\oplus N\ZZ\oplus N\ZZ/4$.
 Then $v(L) = N^3$. Also, $L^* = \ZZ\oplus \ZZ/2\oplus \ZZ/16$ and $B = -32N^3$.

 As a quadratic form, $Q$ is given by the determinant of a matrix, for $x = (x_1,x_2,x_3)$
 \[Q(x) =\ ^txQx = \frac 2N (x_2^2-4x_1x_3) = \frac{-8}N
 \left |\begin{smallmatrix}
 x_1 & x_2/2\\
 x_2/2 & x_3
 \end{smallmatrix}\right |.\]
 And there is a map from $\SL_2(\mathbb{R})$ to $SO(Q)$:
 \[\sm a & b \\ c & d \esm \mapsto
 \sm a^2 & ab & b^2 \\
 2ac & ad + bc & 2bd \\
 c^2 & cd & d^2 \esm.\]
 This map gives an isomorphism of $SO(Q)$ with $\SL_2(\mathbb{R}) / \pm I$.

 \begin{dfn}
 \begin{enumerate}
 \item Let $\Gamma_Q$ be a discrete subgroup of $SO(Q)$ which leaves $L$ invariant.
  Let $\Gamma_Q^*$ be the (normal) subgroup of $\Gamma_Q$ which fixes $L^*/L$.
  Let $\chi$ be a character of $\Gamma_Q$ which is trivial on $\Gamma_Q^*$.
  We say that $\omega: L^*/L\to\CC$ has the {\bf{second permutation property}}
  for $\Gamma_Q$ with character $\chi$ if it satisfies
     \[\omega(\gamma \kappa) = \chi(\gamma)\omega(\kappa),\ \gamma\in\Gamma_Q, \kappa\in L^*.\]

 \item Let $m\in \ZZ$. We say that a function $f\in\mathcal{S}(\mathbb{R}^3)$
 has the {\bf{second spherical property}} for weight $2m$ if it satisfies, by
 identifying $\kappa(\phi)$ as an element of $SO(Q)$
     \[f(\kappa(\phi)^{-1}x) = e^{-2im\phi}f(x)\ \text{for all $\phi\in\mathbb{R}$ and $x\in\mathbb{R}^3$}.\]
 \end{enumerate}
 \end{dfn}

 Let us define Hermite polynomials: for $0\leq \nu\in\ZZ$, define
 \[H_\nu(x) = (-1)^\nu \exp(x^2/2)\frac{d^\nu}{dx^\nu}\exp(-x^2/2).\]

 \begin{thm} \label{spherical} Let $m$ and $\lambda$ be integers.
 Then for every positive integer $\mu$ such that $|m| \leq \lambda + \mu$,
 there is a function $L_{m,\lambda,\mu}$ such that
  \begin{equation}\label{ex1}
  f_{m,\lambda,\mu}(x) = L_{m,\lambda,\mu}(x) H_\mu(\frac{\sqrt{8\pi}}N(x_1+x_3))
  \exp(\frac{-2\pi}N(2x_1^2+x_2^2+2x_3^2))\end{equation}
 has the {\bf{first spherical property}} for weight $k/2 = \lambda+1/2$,
  and the {\bf{second spherical property}} for weight $2m$.
  The function $L_{m,\lambda,\mu}$ is defined by
 \[L_{m,\lambda,\mu}(x) = \frac1{2\pi}\int^{2\pi}_o e^{2mi\phi}L_{\lambda,\mu}(\kappa(\phi)^{-1}x)d\phi\]
 where $L_{\lambda,\mu}(x)
 = H_{\nu_1}(\sqrt{8\pi/N}(x_1-x_3))H_{\nu_2}(\sqrt{8\pi/N}x_2)$ for any choice of $\nu_1$
 and $\nu_2$ such that $\nu_1 + \nu_2 - \mu = \lambda$. In particular, we may take
 \[L_{\lambda,\lambda,0}(x) = (x_1-ix_2-x_3)^\lambda.\]
 \end{thm}

 {\bf{Proof}} This is Theorem 2.1 in \cite{C}. \qed

 \medskip

 Let $f_{k,m} = f_{m,\lambda,\mu}$ for some fixed $\mu$ as in Theorem \ref{spherical}
  where $\frac{k}{2} = \lambda+\frac{1}{2}$. Consider theta function associated
  with $f_{k,m}$: for $z = u+iv ,w = \xi+i\eta \in \HH$ and given character $\chi$
  for $\Gamma_0(4N)$ we define a theta function
 \begin{equation} \label{theta}
 \theta(z,w; f_{k,m}) := (32N^3)^{-1/2}i^\lambda v^{-k/4}
 (4\eta)^{-m}\sum_{x\in L_N^*}\check{\bar{\chi}}_1(4x_1)
 \{r(\sigma_{4Nz},Q)p(\sigma_{2Nw})f_{k,m}\}(x),
 \end{equation}
 where $\chi_1 = \chi(\frac{-1}{})^\lambda$, $\check{\bar{\chi}}_1(l)
 = \displaystyle\sum_{h=1}^{4N}\bar{\chi}_1(h)e^{2\pi hl/4N}$, and $L_N^*
  = \mathbb{Z}/4\oplus \mathbb{Z}/2\oplus \mathbb{Z}/4$,
  the dual lattice to $L_N = N\mathbb{Z}\oplus N\mathbb{Z}\oplus N\mathbb{Z}$.
  \\

 \begin{thm} The above theta function is a two-variable nonholomorphic modular form:
 \begin{enumerate}
 \item $\theta(z,w;f_{k,m})\in G_{k/2}(\Gamma_0(4N),\chi)$ as a function of $z$,
 \item $\overline{\theta(z,w; f_{k,m})} \in G_{2m}(\Gamma_0(2N),\chi^2)$ as a function of $w$.
 \end{enumerate}
 \end{thm}

 {\bf{Proof}} Basically we follow the proof of Theorem 2.3 of \cite{C}.
 Let $z=u+iv, w=\xi+i\eta$ for $u,v,\xi, \eta \in \mathbb{R}$.
 We define
 \[\Theta(z,w; f_{k,m}) = (4\eta)^{-m}v^{-k/4}
 \sum_{x\in L'}\bar{\chi}_1(x_1)\{r(\sigma_z,Q)p(\sigma_{4w})
 f_{k,m}\}(x),\]
 where $L' = \mathbb{Z}\oplus N\mathbb{Z}\oplus N\mathbb{Z}/4$. With the notation of Corollary \ref{first},
 \[\Theta(z,w;f_{k,m})
 = (4\eta)^{-m}\theta(z,p(\sigma_{4w})f_{k,m};\omega)\]
 with $\omega: L^*/L\to\mathbb{C}$ defined by
 \begin{enumerate}
 \item $\omega(\kappa) = 0$ if $\kappa\not\in L'$,
 \item $\omega(\kappa) = \bar{\chi}_1(\kappa_1)$ if $\kappa = (\kappa_1,\kappa_2,\kappa_3)\in L'$,
 \end{enumerate}
 where $L = 4N\mathbb{Z}\oplus N\mathbb{Z} \oplus N\mathbb{Z}/4$
 and its dual is $L^* = \mathbb{Z}\oplus \mathbb{Z}/2\oplus \mathbb{Z}/16$. Note that
 $\omega$ has the first permutation property for $\Gamma_0(4N)$ with character
 $\bar{\chi}_1$, and the second permutation property for $\Gamma_Q=\sm 2&\\&1/2\esm \Gamma_0(2N)\sm 1/2&\\&2\esm$
 with character $\chi^2$ (See Proposition 2.2 in \cite{C}). Note that $p(\sigma_{4w})f_{k,m}$ has the first
  spherical property of weight $k/2$ since $p$ commutes with the Weil representation (See (\ref{commute})).
  Then by the Corollary \ref{first}, $\Theta(z,w; f_{k,m}) \in G_{k/2}(\Gamma_0(4N), \bar{\chi}(\frac N{}))$
 as a function of $z$.
\\

 If we use equation (\ref{p}) then we see that for $\gamma = \sm a&b\\c&d\esm \in\Gamma_0(2N)$
 \[p(\sigma_{4(\gamma w)})(x) = p(\gamma_2\sigma_{4w}\kappa(\phi)^{-1})(x)\]
 where $\gamma_2 = \sm 2&\\&1/2\esm \gamma \sm 1/2&\\&2\esm$ and $e^{-i\phi} = J(\gamma,w)/|J(\gamma,w)|$.
 Note that $\gamma_2\in \Gamma_Q$. So if we use the second spherical property of $f_{k,m}$ and second
 permutation property of $\omega$ then we get the desired transformation property for the second variable $w$.
\\

 We define the Fricke involution $W(N)$:
 \begin{equation*}
 (f|_k W(N))(z)=
 \begin{cases}
 N^{-k/2}(-iz)^{-k}f(-1/Nz)& \text{$k=$ half-integer},\\
 N^{-k/2}z^{-k}f(-1/Nz) & \text{$k=$ even-integer}.
 \end{cases}
 \end{equation*}
 Let $|_{k/2}W(4N)$ act on the variable $z$, and $|_{2m}W(2N)$ act on $w$ (See \cite{C}).
 Then
 \[\theta(z,w;f_{k,m}) = (\Theta|_{k/2} W(4N) \overline{|_{2m} W(2N)})(z,w;f_{k,m})\]
 and hence transformation formulas of $\theta(z,w;f_{k,m})$ follow. \qed
\\

 \medskip

 This theta function also has a good property about differential operators.

 \begin{prop} \label{PDE} The theta function $\theta(z,w;f_{k,m})$,
 defined in (\ref{theta}), satisfies the PDE
 \[4[v^2(\frac{\partial^2}{\partial u^2}+\frac{\partial^2}{\partial v^2})
 -i\frac k2 v(\frac \partial{\partial u}+i\frac \partial{\partial v})+\frac k4(\frac k4 -1)]\theta(z,w;f_{k,m})\]
 \[=[\eta^2(\frac{\partial^2}{\partial\xi^2}+\frac{\partial^2}{\partial\eta^2})
 +2mi\eta(\frac \partial{\partial\xi}-i\frac \partial{\partial\eta})+m(m-1) - \frac 34]\theta(z,w;f_{k,m}).\]
 \end{prop}
 {\bf{Proof}} This is Proposition 2.13 of \cite{C}. \qed

 \subsection{{\bf{Regularized theta lifting}}}

 In this section we explain how to regularize the integral and we define theta lifting using that regularization.
  We will use Borcherds' {\bf{regularized integral}}. Throughout we use the setup of \cite{B1}. If we use the
  weak Maass form, the integral in the theta lifting may
   be divergent. So it has to be regularized as follows:
   We integrate over the region $\FF_t$, where $\FF_\infty = \{z\in\HH |\ |z|\geq1, |\re(z)|\leq 1/2\}$ is
    the usual fundamental domain of $\SL_2(\mathbb{Z})$ and $\FF_t$ is the subset of $\FF_\infty$
    of points $z$ with $\im(z)\leq t$. Suppose that
 \begin{equation} \label{limit}
 \lim_{t\rightarrow \infty}\int_{\FF_t}F(z)v^{-s}\frac{dudv}{v^2}
 \end{equation}
 exists for $\re(s)\gg0$ and can be continued to a meromorphic
 function defined for all complex $s$. Then we define
 \[\int^{\mathrm{reg}}_{\SL_2(\ZZ)\setminus\HH}F(z)\frac{dudv}{v^2}\]
 to be the constant term of the Laurent expansion of the function (\ref{limit}) at $s=0$. If we use this
  regularized integral, we can define the theta lifting even though we use weak Maass forms.
  \\

 This is a regularized integral for $\SL_2(\ZZ)$ but we need a regularized integral for $\Gamma_0(4N)$.
 Note that for a modular form $g\in S_{k/2}(\Gamma_0(4N),\chi)$
 \[\int_{\Gamma_0(4N)\setminus\HH}v^{k/2}g(z)\overline{\theta(z,w;f_{k,m})}\frac{dudv}{v^2} = \int_{\SL_2(\ZZ)
 \setminus\HH}v^{k/2}\sum_{\alpha\in\Gamma_0(4N)\setminus\SL_2(\ZZ)}
 g_\alpha(z)\overline{\theta_\alpha(z,w;f_{k,m})}\frac{dudv}{v^2}\]
 where $g_\alpha(z) = (cz+d)^{-k/2}g(\alpha z)$ for $\alpha = \sm a&b\\c&d\esm$. For $g\in
 WMF_{k/2,s}^*(\Gamma_0(4N),\chi)$ define a function on $\HH$
 \[\Phi(g)(w) = \int^{\mathrm{reg}}_{\SL_2(\ZZ)\setminus\HH}v^{k/2}
 \sum_{\alpha\in\Gamma_0(4N)\setminus\SL_2(\ZZ)}
 g_\alpha(z)\overline{\theta_\alpha(z,w;f_{k,m})}\frac{dudv}{v^2}.\]
\\

 And for square-free, positive integer $D$ and $g\in WMF^*_{k/2,s}(\Gamma_0(4N),\chi)$ let $g_D(x) = g(Dx)$.
 Then $g_D\in WMF^*(\Gamma_0(4ND),\chi_D)$. Define
 \[\Phi_D(g)(w) = \Phi(g_D)(\frac wD).\]
 But in this case, we use $\theta(z,w;f_{k,m})$ such that its level is $4ND$ and its character is
 $\chi_D = \chi(\frac{-1}{})^\lambda(\frac{D}{})$ and use a regularized integral for $\Gamma_0(4ND)$.
\bigskip

\section{{\bf{Proof of Theorem\ref{main1}}}}

 First we prove the convergence of $\Phi_D(g)$ for $g\in
 WMF^*_{k/2,s}(\Gamma_0(4N),\chi)$. Since $\Phi_D(g)$
 is defined by using $\Phi(g_D)$, it is enough to show the
 convergence of $\Phi(g)$. Note that since
 $\displaystyle\sum_{\alpha\in\Gamma_0(4N)\setminus\SL_2(\ZZ)}$
 is a finite sum, we can exchange sum
 and integration. And since $\int_{\FF_t} = \int_{\FF_1} +
 \int_{v=1}^{v=t}\int_{u=-1/2}^{u=1/2}$, we only need to check
 \begin{equation} \label{convergence}
 \lim_{t\to\infty}\int_{v=1}^{v=t}\int_{u=-1/2}^{u=1/2}v^{k/2}
 g_{\alpha}(z)\overline{\theta_\alpha(z,w;f_{k,m})}\frac{dudv}{v^2}
 \end{equation}
 for each $\alpha\in\Gamma_0(4N)\setminus\SL_2(\ZZ)$. And if
 $g_\alpha(z)\overline{\theta_\alpha(z,w;f_{k,m})}$
 has a Fourier expansion as $\displaystyle\sum_{n\in\ZZ}a(n,v,w)e^{2\pi inu}$ then
 \[\int_{v=1}^{v=t}\int_{u=-1/2}^{u=1/2}v^{k/2}g_{\alpha}(z)
 \overline{\theta_\alpha(z,w;f_{k,m})}\frac{dudv}{v^2}
 = \int_{v=1}^{v=t}v^{k/2-2}a(0,v,w)dv.\]
 So we need to check the Fourier coefficients of the constant term
 of $g_\alpha(z)\overline{\theta_\alpha(z,w;f_{k,m})}$.
 For this we need Fourier expansions of $g_\alpha(z)$ and $\theta_\alpha(z,w;f_{k,m})$ with respect to $z$.

 By computing explicitly representations in the definition of $\theta(z,w;f_{k,m})$ in (\ref{theta}) it turns out that
 \[\theta(z,w;f_{k,m}) = \sum_{x\in\ZZ^3} h(x,v,w;k)e^{\frac{-\pi v}{4N^2}|\Lambda(x,w)|^2}e^{2\pi i\bar{z}(x_2^2-x_1x_3)}\]
 where $h$ is a polynomial of $x,v$ and $w$ and $\Lambda(x,w) = \frac{1}{4\eta}(x_1-4Nwx_2+4N^2w^2x_3)$.
 This gives a Fourier expansion with respect to $z=u+iv$.

 Next we will see the Fourier expansion of $\theta_\alpha(z,w;f_{k,m})$
 for general $\alpha\in\Gamma_0(4N)\setminus\SL_2(\ZZ)$. Let
 \[Q_4 = \frac12 \sm &&-2\\&1&\\-2&&\esm.\]
 This is just the original $Q$ in (\ref{Q}) when $N=4$. Likewise let $f_4$ be a function $f_{k,m}$ in (\ref{ex1})
  with $N=4$. Then $f_4$ satisfies the first and second spherical properties for the weights $k/2$ and $2m$ respectively. Now let
 \begin{eqnarray*}
 L &=& 4N\ZZ \oplus 2\ZZ\oplus \ZZ\\
 L' &=& \ZZ \oplus 2\ZZ \oplus \ZZ\\
 L^* &=& \ZZ\oplus \ZZ\oplus \ZZ/4N
 \end{eqnarray*}
 and define $\omega: L^*/L\to\CC$ by
 \begin{equation*}
 \omega(\kappa) =
 \begin{cases}
 0 & \kappa\not\in L'\\
 \check{\bar{\chi}}(\kappa_1) & \kappa = (\kappa_1,\kappa_2,\kappa_3) \in L'.
 \end{cases}
 \end{equation*}
 Then $\omega$ has the first permutation property for $\Gamma_0(4N)$ with character $\chi_1$ and the second
 permutation property for $\Gamma_Q=\sm 2&\\&1/2\esm \Gamma_0(2N)\sm 1/2&\\&2\esm$ with character $\bar{\chi}^2$.
 Direct computation shows that, up to constant multiple,
 \[\theta(z,w;f_{k,m}) \approx \eta^{-m}\theta(z,p(\sigma_{2Nw})f_4;\omega).\]
 By Theorem \ref{Shintani} we see that for each $\alpha\in\Gamma_0(4N)\setminus\SL_2(\ZZ)$,
  $\theta_\alpha(z,w;f_{k,m})$ can be written as
 \[\theta_\alpha(z,w;f_{k,m}) =
 \sum_{x\in\ZZ^3}h_\alpha(x,v,w;k)e^{\frac{-4\pi v}{4N^2}|\Lambda(x,w)|^2}e^{2\pi i\bar{z}(x_2^2-x_1x_3)}\]
 where $h_\alpha$ is another polynomial for $x,v$ and $w$.

 And we know that $g_\alpha$ has the Fourier expansion of the form (\ref{fourier0}). So the constant term
 of $g_\alpha(z)\overline{\theta_\alpha(z,w;f_{k,m})}$ is a sum of terms of the form $a(x_2^2-x_1x_3,v)
 e^{2\pi(x_2^2-x_1x_3)v}h_\alpha(x,v,w;k)e^{\frac{-\pi v}{4N^2} | \Lambda(x,w) |^2}$ where $x=(x_1,x_2,x_3)
 \in\mathbb{Z}^3$. By the definition of $WMF^*_{k,s}(\Gamma_0(4N),\chi)$, $a(x_2^2-x_1x_3,v)e^{2\pi(x_2^2-x_1x_3)v}
 = O(v^\delta)$ for some $\delta>0$. So every term is exponentially
 decreasing as $v\to\infty$ except the case of
 $\Lambda(x,w) =0$. So if $\Lambda(x,w)\neq 0$ for all $x\in\ZZ^3$
 then the constant term of $g_\alpha(z)\overline{\theta_\alpha(z,w;f_{k,m})}$
 goes to zero when $v$ goes to $\infty$. So (\ref{convergence})
 converges and hence $\Phi(g)(w)$ is well defined where $\Lambda(x,w)\neq 0$
 for all $x\in\ZZ^3$. And it may have singularities where $\Lambda(x,w) = 0$
 for some $x\in\ZZ^3$. Since $\chi(0)=0$, we don't need to
 consider the case of $x=0$. So the singularities may occur where
 $w = \frac{b\pm\sqrt{b^2-ac}}{2Nc}\in\HH$ for $a,b,c\in\ZZ$.

 And transformation properties of $\Phi(g)(w)$ come easily from the fact that $\overline{\theta(z,w;f_{k,m})}$
 is a nonholomorphic modular form of weight $2m$ for $\Gamma_0(2N)$ and $\chi^2$ as a function of $w$.

 We have Maass differential operators on smooth functions from $\mathbb{H}$ to $\mathbb{C}$ (See \cite{B3} page 97):
 \begin{eqnarray*}
 R_k &=& 2i\frac\partial{\partial z} + kv^{-1},\\
 L_k &=& 2iv^2\frac\partial{\partial\bar{z}}.
 \end{eqnarray*}
 For any smooth function $f: \mathbb{H} \rightarrow \mathbb{C}$ and $\gamma\in \SL_2(\mathbb{Z})$ it is well known that
 \begin{eqnarray*}
 (R_k f)|_{k+2}\gamma &=& R_k(f|_{k}\gamma),\\
 (L_k f)|_{k-2}\gamma &=& L_k(f|_{k}\gamma).
 \end{eqnarray*}
 The operator $\Delta_k$ can be expressed in terms of $R_k$ and $L_k$ by
 \[\Delta_k = L_{k+2}R_k - k = R_{k-2}L_k.\]

 We will use $\Delta_{2m}$ as the Laplace operator
 with respect to $w$ and $\Delta_{k/2}$ as the Laplace operator with respect to $z$.
 To prove that $\Phi(g)(w)$ is an eigenfunction of $\Delta_{2m}$
  we need following lemmas, which are essentially Lemma 4.2 and Lemma 4.3 in \cite{B3}.

 \begin{lem} \label{previous} Let $f\in G_{k/2}(\Gamma_0(4N),\chi)$ and
 $g\in G_{k/2+2}(\Gamma_0(4N),\chi)$. Then
 \begin{eqnarray*}
 &&\int_{\mathcal{F}_t}\sum_{\alpha\in \Gamma_0(4N)\setminus \SL_2(\mathbb{Z})}
 f_\alpha(z)\overline{(L_{k/2+2}(g_\alpha))(z)}v^{k/2-2}dudv\\
 &-& \int_{\mathcal{F}_t}\sum_{\alpha\in \Gamma_0(4N)\setminus
 \SL_2(\mathbb{Z})}(R_{k/2}(f_\alpha))(z)\overline{g_\alpha(z)}v^{k/2} dudv\\
 &=& \int^{\frac12}_{-\frac12}[\sum_{\alpha\in \Gamma_0(4N)\setminus \SL_2(\mathbb{Z})}
 f_\alpha(z)\overline{g_\alpha(z)}v^{k/2}]_{v=t}du.
 \end{eqnarray*}
 \end{lem}

 {\bf{Proof}} The assumptions imply that $\displaystyle\omega = v^{k/2}
 \big(\sum_{\alpha\in \Gamma_0(4N)\setminus \SL_2(\mathbb{Z})}
 f_\alpha(z)\overline{g_\alpha(z)}\big)d\bar{z}$ is a $\SL_2(\mathbb{Z})$
 -invariant $1$-form on $\mathbb{H}$. By Stokes' theorem we have
 \begin{eqnarray*}
 &&\int_{\partial\mathcal{F}_t}v^{k/2}\sum_{\alpha\in \Gamma_0(4N)\setminus
 \SL_2(\mathbb{Z})}f_\alpha(z)\overline{g_\alpha(z)}d\bar{z}\\
 &=& \int_{\mathcal{F}_t}d(v^{k/2}\sum_{\alpha\in \Gamma_0(4N)\setminus
 \SL_2(\mathbb{Z})}f_\alpha(z)\overline{g_\alpha(z)}d\bar{z})\\
 &=& \int_{\mathcal{F}_t}\big(-\frac{\partial}{\partial v}v^{k/2}
 \sum_{\alpha\in \Gamma_0(4N)\setminus \SL_2(\mathbb{Z})}f_\alpha(z)
 \overline{g_\alpha(z)}-i\frac{\partial}{\partial u}v^{k/2}
 \sum_{\alpha\in \Gamma_0(4N)\setminus \SL_2(\mathbb{Z})}f_\alpha(z)
 \overline{g_\alpha(z)}\big)dudv\\
 &=&\int_{\mathcal{F}_t}\big(v^{k/2-2}\sum_{\alpha\in \Gamma_0(4N)\setminus
  \SL_2(\mathbb{Z})}f_\alpha(z)\overline{(L_{k/2+2}(g_\alpha))(z)}\\
 &&-v^{k/2}\sum_{\alpha\in \Gamma_0(4N)\setminus \SL_2(\mathbb{Z})}
 (R_{k/2}(f_\alpha))(z)\overline{g_\alpha(z)}\big)dudv.
 \end{eqnarray*}
 In the integrand over $\partial\mathcal{F}_t$ on the left hand side
 the contributions from $\SL_2(\mathbb{Z})$-equivalent boundary pieces cancel. Thus
 \[\int_{\partial\mathcal{F}_t}v^{k/2}\sum_{\alpha\in \Gamma_0(4N)\setminus
  \SL_2(\mathbb{Z})}f_\alpha(z)
 \overline{g_\alpha(z)}d\bar{z} = \int^{\frac12}_{-\frac12}[v^{k/2}
 \sum_{\alpha\in \Gamma_0(4N)\setminus \SL_2(\mathbb{Z})}f_\alpha(z)\overline{g_\alpha(z)}]_{v=t}du.\]
 This implies the assertion. \qed

 \begin{lem} \label{Bruinier} Let $f,g\in G_{k/2}(\Gamma_0(4N),\chi)$. Then
 \begin{eqnarray*}
 &&\int_{F_t}v^{k/2}\sum_{\alpha\in \Gamma_0(4N)\setminus \SL_2(\mathbb{Z})}
 (\Delta_{k/2}(f_\alpha))(z)\overline{g_\alpha(z)}\frac{dudv}{v^2}\\
 &-&\int_{F_t}v^{k/2}\sum_{\alpha\in \Gamma_0(4N)\setminus \SL_2(\mathbb{Z})}
 f_\alpha(z)\overline{(\Delta_{k/2}(g_\alpha))(z)}\frac{dudv}{v^2}\\
 &=&\int^{1/2}_{-1/2}[v^{k/2}\sum_{\alpha\in \Gamma_0(4N)\setminus
 \SL_2(\mathbb{Z})}f_\alpha(z)\overline{(L_{k/2}(g_\alpha))(z)}]_{v=t}du\\
 &-&\int^{1/2}_{-1/2}[v^{k/2}\sum_{\alpha\in \Gamma_0(4N)\setminus
 \SL_2(\mathbb{Z})}(L_{k/2}(f_\alpha))(z)\overline{g_\alpha(z)}]_{v=t}du.
 \end{eqnarray*}
 \end{lem}

 {\bf{Proof}} We write $\Delta_{k/2} = R_{k/2-2}L_{k/2}$ and apply Lemma \ref{previous} twice. {\qed}

 \medskip

 Observe that $\Delta_{k/2}(g_\alpha) = sg_\alpha$. We have
 \begin{eqnarray*}
 \Delta_{2m}(\Phi(g))(w)
 &=& -4\int^{\mathrm{reg}}_{\SL_2(\ZZ)\setminus\HH}v^{k/2}\sum_{\alpha\in
 \Gamma_0(4N)\setminus \SL_2(\mathbb{Z})}g_\alpha(z)
 \overline{(\eta^2\frac{\partial^2}{\partial w\partial\bar{w}}+2m i
 \eta\frac{\partial}{\partial w})\theta_\alpha(z,w;f_{k,m})}\frac{dudv}{v^2}\\
 &=& 4\int^{\mathrm{reg}}_{\SL_2(\ZZ)\setminus\HH}v^{k/2}
 \sum_\alpha g_\alpha(z)\overline{(\Delta_{k/2}(\theta_\alpha))(z,w;f_{k,m})}\frac{dudv}{v^2}\\
  &&+
  4\big(m(m-1)-\frac34-k(\frac k4-1)\big)\Phi(g)(w).
 \end{eqnarray*}
 The last equality comes from Proposition \ref{PDE}. And by Lemma \ref{Bruinier},
 \begin{eqnarray*}
 && \int^{\reg}_{\SL_2(\ZZ)\setminus\HH}v^{k/2}
 \sum_\alpha g_\alpha(z)
 \overline{(\Delta_{k/2}(\theta_\alpha))(z,w;f_{k,m})}\frac{dudv}{v^2}\\
 &=& s\Phi(g)(w)-\lim_{v\rightarrow\infty}\int^{1/2}_{-1/2}v^{k/2}\sum_\alpha
 g_\alpha(z)\overline{(L_{k/2}(\theta_\alpha))(z,w;f_{k,m})}du\\
 && +\lim_{v\rightarrow\infty}\int^{1/2}_{-1/2}v^{k/2}
 \sum_\alpha(L_{k/2}(g_\alpha))(z)\overline{\theta_\alpha(z,w;f_{k,m})}du.
 \end{eqnarray*}

 By the same argument in the proof of the convergence, we see that two integrals
 \begin{eqnarray*}
 &&\lim_{v\rightarrow\infty}\int^{1/2}_{-1/2}v^{k/2} \sum_\alpha g_\alpha(z)\overline{(L_{k/2}\theta_\alpha)(z,w;f_{k,m})}du,\\
 &&\lim_{v\rightarrow\infty}\int^{1/2}_{-1/2}v^{k/2} \sum_\alpha (L_{k/2}(g_\alpha))(z)\overline{\theta_\alpha(z,w;f_{k,m})}du
 \end{eqnarray*}
 vanish.

 And note that $\Phi(g_D)(w)$ transforms under the group $\Gamma_0(2ND)$. So $\Phi(g_D)(w/D)$
 transforms under the group $\Gamma_0(2N,D)$.
 \bigskip

 \section{\bf{Proof of Theorem \ref{main2}}}

 Let $\lambda = (k-1)/2$. We consider the lifting from $H^+_{k/2}(\Gamma_0(4N),\chi)$
 to $H^*_{2\lambda}(\Gamma_0(2N,D),\chi^2)$
 with singularities. In this case we take $f_{k,\lambda}(x) = (x_1-ix_2-x_3)^\lambda
  \exp(-\frac{2\pi}{N}(2x_1^2+x_2^2+2x_3^2))$.
  Then $\theta(z,w;f_{k,\lambda})$ is the same theta function which was used by Niwa\cite{N}
  and Cipra\cite{C}.
\\

 As a special case of Theorem \ref{main1}, if $g\in H^+_{k/2}(\Gamma_0(4N),\chi)$
 then $\Phi_D(g)(w)\in H_{2\lambda}(\Gamma_0(2N,D),\chi^2)$
 with singulairites. Now we will compute the image of the lifting when $k\geq1$
 by using the unfolding method. To do that we need to rewrite theta functions as follows:

 \begin{lem} \label{unfolding} We have
 \begin{eqnarray*}
 \theta(z,i\eta;f_{k,\lambda}) &=& C\sum_{\nu = 0}^\lambda \sm \lambda \\ \nu\esm
 (2/\pi)^{\nu/2}\eta^{1-\nu}\sum_{\gamma\in\Gamma_\infty\setminus\Gamma_0(4N)}
 \frac{\chi(d)}{(\im\gamma z)^{\lambda-\nu/2}j(\gamma,z)^k}\\
 &&\times \sum_{m,n=-\infty}^\infty \bar{\chi}_1(m)m^{\lambda-\nu}H_\nu(2\sqrt{2\pi
 \im\gamma z}n)e^{2\pi in^2\gamma z-\frac{\pi m^2\eta^2}{4\im \gamma z}}
 \end{eqnarray*}
 with $C = (-1)^\lambda 2^{-4\lambda}N^{\lambda/2-1/4}$.
 \end{lem}

 {\bf{Proof}} This is Theorem 2.11 of \cite{C}. \qed

 \medskip

 Now we will observe the behavior of $\Phi(g)(w)$ at $i\infty$. The following computations
 are essentially done in Theorem 2.12 in \cite{C}. Assume that $k\geq1$. Let
 \[g(z) = \sum_{n\gg -\infty}a^+(n)e(nz) + \sum_{n<0}a^-(n)W_{k/2}(2\pi nv)e(nz)\in H^+_{k/2}(\Gamma_0(4N),\chi).\]
 Using the Lemma \ref{unfolding}, we have
 \begin{eqnarray*}
 \Phi(g)(i\eta) &=& \int_{\Gamma_0(4N)\setminus\mathbb{H}}v^{k/2}g(z)
 \overline{\theta(z,i\eta;f_{k,\lambda})}\frac{dudv}{v^2}\\
 &=&C\sum_{\nu=0}^\lambda \sm \lambda\\ \nu\esm (2/\pi)^{\nu/2}\eta^{1-\nu}
 \int^\infty_0\int^1_0 v^{k/2}g(z)v^{\nu/2-\lambda}\sum_{m,n=-\infty}^\infty \chi_1(m)\\
 &&\times m^{\lambda-\nu}H_\nu(2\sqrt{2\pi v}n)e^{-2\pi in^2\bar{z}-
 \frac{\pi\eta^2m^2}{4v}}\frac{dudv}{v^2}\\
 &=& C\sum_{\nu=0}^\lambda \sm \lambda\\ \nu\esm (2/\pi)^{\nu/2}\eta^{1-\nu}
 \int^\infty_0v^{(\nu-1)/2}\sum_{n=-\infty}^\infty a^+(n^2)H_\nu(2\sqrt{2\pi xv}n)\\
 &&\times e^{-4\pi n^2v}\sum_{m=-\infty}^\infty \chi_1(m)m^{\lambda-\nu}e^{-\pi\eta^2m^2/4v}dv/v\\
 &=&C'\sum_{\nu=0}^\lambda\sm \lambda\\ \nu\esm (2\pi)^{-\nu}\int^\infty_0[a^+(0)H_\nu(0)+\sum_{n\neq0}a^+(n^2)H_\nu(yn)e^{-n^2y^2/2}]\\
 &&\times (\eta/y)^{1-\nu}\sum_{m=-\infty}^\infty \chi_1(m)m^{\lambda-\nu}e^{-2\pi^2m^2(\eta/y)^2}dy/y
 \end{eqnarray*}
 with $C' = 2C(8\pi)^{1/2} = (-1)^\lambda2^{-4\lambda+2}N^{\lambda/2+1/4}(2\pi)^{1/2}$. As $\eta\to\infty$,
 the only non-negligible term is that one involving $a^+(0)$:
 \begin{eqnarray*}
 \Phi(g)(i\infty) &=& a^+(0)C'\sum_{\nu=0}^\lambda\sm \lambda\\ \nu\esm(2\pi)^{-\nu}H_\nu(0)\int^\infty_0y^{1-\nu}
 \sum_{m=-\infty}^\infty \chi_1(m)m^{\lambda-\nu}e^{-2\pi^2m^2y^2}dy/y.
 \end{eqnarray*}
 If $k\geq3$ then we invert the theta function using Poisson summation
 \begin{eqnarray*}
 \Phi(g)(i\infty) &=& a^+(0)C''\sum_{\nu=0}^\lambda\sm \lambda\\ \nu\esm i^\nu H_\nu(0)
 \int^\infty_0y^{-\lambda}\sum_{m=-\infty}^\infty\check{\chi}_1(m)H_{\lambda-\nu}(\frac m{4Ny})e^{-m^2/32N^2y^2}dy/y
 \end{eqnarray*}
 with $C'' = C'(2\pi i)^{-\lambda}(2\pi)^{-1/2}(4N)^{-1} = i^\lambda 2^{-5\lambda}
 N^{\lambda/2-3/4}\pi^{-\lambda}$. We can now actually sum over $\nu$:
 \[\sum_{\nu=0}^\lambda\sm \lambda\\ \nu\esm i^\nu H_\nu(0)H_{\lambda-\nu}(\frac m{4Ny}) = (\frac m{4Ny})^\lambda\]
 so
 \begin{eqnarray*}
 \Phi(g)(i\infty) &=& a^+(0)C''\int^\infty_0\sum_{m=-\infty}^\infty\check{\chi}_1(m)(\frac m{4Ny^2})^\lambda e^{-m^2/32N^2y^2}dy/y\\
 &=& a^+(0)C^*\int^\infty_0\sum_{m=1}^\infty\frac{\check{\chi}_1(m)}{m^\lambda}y^\lambda e^{-\nu}dy/y
 \end{eqnarray*}
 with $C^* = C''(8N)^\lambda = i^\lambda2^{-3\lambda}N^{(3/2)\lambda-(3/4)}\pi^{-\lambda}$.
 The integral gives $\Gamma(\lambda)$, and $C^*\Gamma(\lambda) = C_0(\lambda)$. So
 \[\Phi(g)(i\infty) = a^+(0)C_0(\lambda)\sum_{m=1}^\infty \check{\chi}_1(m)m^{-\lambda}\]
 where $C_0(\lambda) = i^\lambda 2^{-2\lambda}N^{(3/2)\lambda-(3/4)}\pi^{-\lambda}\Gamma(\lambda)$
 and $\check{\chi}_1(m) = \sum_{h=1}^{4N}\chi_1(h)e^{\pi imh/2N}$.
   And if we use functional equations for the $L$-series, we get the result about the constant term.

 If $k=1$ then we can evaluate directly without the inversion
 \begin{eqnarray*}
 \Phi(g)(i\infty) &=& a^+(0)C'\int^\infty_0 y\sum_{-\infty}^\infty \chi(m)\mathrm{exp}(-2\pi^2m^2y^2)dy/y\\
 &=& a^+(0)C'2^{-1/2}\pi^{-1}\Gamma(1/2)\sum_{m=1}^\infty\frac{\chi(m)}{m}\\
 &=&4N^{1/4}a^+(0)\sum_{m=1}^\infty \frac{\chi(m)}{m}.
 \end{eqnarray*}

 From this, we see that $\Phi(g)(w)$ is bounded at $i\infty$. So we see that $\Phi(g)(w)\in
 H^*_{2\lambda}(\Gamma_0(2N),\chi^2)$ with singularities and $P(\Phi(g))$ is constant. We have
 \begin{eqnarray*}
 \Phi(g)(i\eta)-\Phi(g)(i\infty) &=& C'\sum_{\nu=0}^\lambda\sm \lambda\\ \nu
 \esm (2\pi)^{-\nu}\int^\infty_0\sum_{n\neq0}a^+(n^2)H_\nu(yn)e^{-n^2y^2/2}(\eta/y)^{1-\nu}\\
 &&\times \sum_{m=-\infty}^\infty \chi_1(m)m^{\lambda-\nu}e^{-2\pi^2m^2(\eta/y)^2}dy/y.
 \end{eqnarray*}
 Now
 \begin{eqnarray*}
 &&\int^\infty_0 \eta^{s-1}(\eta/y)^{1-\nu}\sum_{m=-\infty}^\infty\chi_1(m)m^{\lambda-\nu}e^{-2\pi^2m^2(\eta/y)^2}d\eta\\
 &=&y^s(2\pi^2)^{(\nu-s-1)/2}\sum_{m=1}^\infty\chi_1(m)m^{-(s-\lambda+1)}\int^\infty_0\mu^{(s+1-\nu)/2}e^{-\mu}d\mu/\mu\\
 &=&y^s(2\pi^2)^{(\nu-s-1)/2}L(s-\lambda+1,\chi_1)\Gamma(\frac{s+1-\nu}2).
 \end{eqnarray*}
 Thus
 \[\int^\infty_0 \eta^{s-1}(\Phi(g)(i\eta)-\Phi(g)(i\infty))d\eta = C'L(s-\lambda+1,\chi_1)
 \sum_{\nu=0}^\lambda\sm \lambda\\ \nu\esm(2\pi)^{-\nu}(2\pi^2)^{(\nu-s-1)/2}\]
 \[\times\Gamma(\frac{s+1-\nu}2)\int^\infty_0 y^s\sum_{n\neq0}a^+(n^2)H_\nu(yn)e^{-n^2y^2/2}dy/y.\]
 Now
 \begin{eqnarray*}
 &&\int^\infty_0 y^s\sum_{n\neq0}a^+(n^2)H_\nu(yn)e^{-n^2y^2/2}dy/y\\
 &=&(\sum_{n=1}^\infty \frac{a^+(n^2)}{n^s})\int^\infty_0 y^s(H_\nu(y)+H_\nu(-y))e^{-\nu^2/2}dy/y\\
 &=&
 \begin{cases}
 0 & \text{if $\nu$ is odd}\\
 \displaystyle2(\sum_{n=1}^\infty\frac{a^+(n^2)}{n^s})\int^\infty_0y^{s-1}(-1)^\nu(\frac{d^\nu}{dy^\nu}e^{-y^2/2})dy & \text{if $\nu$ is even}
 \end{cases}\\
 &=&
 \begin{cases}
 0\\
 \displaystyle2(\sum^\infty_{n=1}\frac{a^+(n^2)}{n^s})(s-1)\cdots(s-\nu)\int^\infty_0y^{s-\nu-1}e^{-y^2/2}dy
 \end{cases}\\
 &=&
 \begin{cases}
 0\\
 \displaystyle(\sum^\infty_{n=1}\frac{a^+(n^2)}{n^s})(s-1)\cdots(s-\nu)\Gamma(\frac{s-\nu}2)2^{(2-\nu)/2}.
 \end{cases}
 \end{eqnarray*}

 Thus
 \[\int^\infty_0\eta^{s-1}(\Phi(g)(i\eta)-\Phi(g)(i\infty))d\eta = C'L(s-\lambda+1,\chi_1)
 (\sum^\infty_{n=1}\frac{a^+(n^2)}{n^s})\sum_{\nu\ \text{even}}\sm\lambda\\ \nu\esm\]
 \[\times\frac{\Gamma(\frac{s+1-\nu}2)\Gamma(\frac{s-\nu}2)}{2^{\nu+1/2}\pi^{s+1}}(s-1)\cdots(s-\nu).\]
 Using the identity $\Gamma(t/2)\Gamma((t+1)/2) = 2^{1-t}\pi^{1/2}\Gamma(t)$, we get
 \[\int^\infty_0\eta^{s-1}(\Phi(g)(i\eta)-\Phi(g)(i\infty))d\eta = C''L(s-\lambda+1,\chi_1)
 (\sum^\infty_{n=1}\frac{a^+(n^2)}{n^s})\sum_{\nu\ \text{even}}\sm\lambda\\ \nu\esm \]
 \[\times\frac{\Gamma(s-\nu)(s-1)\cdots(s-\nu)}{(2\pi)^s}\]
 with $C'' = 2C'(2\pi)^{-1/2} = (-1)^\lambda2^{-4\lambda+2}N^{\lambda/2+1/4}$. Note that
  $\Gamma(s-\nu)(s-1)\cdots(s-\nu) = \Gamma(s)$ and $\displaystyle\sum_{\nu\ \text{even}}\sm\lambda\\ \nu\esm = 2^{\lambda-1}$. So
 \[\int^\infty_0 \eta^{s-1}(\Phi(g)(i\eta)-\Phi(g)(i\infty))d\eta = C_1(\lambda)(2\pi)^{-s}
 \Gamma(s)L(s-\lambda+1,\chi_1)\sum_{n=1}^\infty\frac{a^+(n^2)}{n^s}\]
 with $C_1(\lambda) = (-1)^\lambda 2^{-3\lambda+2}N^{\lambda/2+1/4}$.

\end{document}